\documentclass[a4paper,12pt]{amsart}

\usepackage{amsmath}
\usepackage{amssymb}
\usepackage{mathrsfs}
\usepackage{ifthen}
\usepackage{graphicx}
\usepackage[T1]{fontenc} 

\nonstopmode \numberwithin{equation}{section}
\newtheorem{thm}{Theorem}[section]
\newtheorem{cor}[equation]{Corollary}
\newtheorem{lem}[equation]{Lemma}

\newtheorem{Cor}{Corollary}

\theoremstyle{definition}

\newtheorem{prob}[equation]{Problem}
\newtheorem{rem}{Remark}[section]

\newcounter{minutes}
\divide\time by 60
\newcounter{hours}
\multiply\time by 60
\addtocounter{minutes}{-\time}

\newcounter {own}
\def\theown {\thesection       .\arabic{own}}

\newenvironment{pf}[1][]{%
 \vskip 3mm
 \noindent
 \ifthenelse{\equal{#1}{}}%
  {{\slshape Proof. }}%
  {{\slshape #1.} }%
 }%
{\qed\bigskip}

\newcounter{alphabet}

\newenvironment{Thm}[1][]{\refstepcounter{alphabet}%
\bigskip%
\noindent%
{\bf Theorem \Alph{alphabet}}%
\ifthenelse{\equal{#1}{}}{}{ (#1)}%
{\bf .} \itshape}{\vskip 8pt}

\def\be{\begin{equation}}
\def\ee{\end{equation}}

\newcommand{\bee}{\begin{enumerate}}
\newcommand{\eee}{\end{enumerate}}

\newcommand{\blem}{\begin{lem}}
\newcommand{\elem}{\end{lem}}
\newcommand{\bthm}{\begin{thm}}
\newcommand{\ethm}{\end{thm}}
\newcommand{\bcor}{\begin{cor}}
\newcommand{\ecor}{\end{cor}}
\newcommand{\beg}{\begin{examp}}
\newcommand{\eeg}{\end{examp}}
\newcommand{\begs}{\begin{examples}}
\newcommand{\eegs}{\end{examples}}
\newcommand{\bdefe}{\begin{defin}}
\newcommand{\edefe}{\end{defin}}
\newcommand{\bprob}{\begin{prob}}
\newcommand{\eprob}{\end{prob}}
\newcommand{\bei}{\begin{itemize}}
\newcommand{\eei}{\end{itemize}}

\begin{document}

\title{On the Generalized Zalcman Conjecture}

\author{Vasudevarao Allu}
\address{Vasudevarao Allu,
	Discipline of Mathematics,
	School of Basic Sciences,
	Indian Institute of Technology  Bhubaneswar,
	Argul, Bhubaneswar, PIN-752050, Odisha (State),  India.}
\email{avrao@iitbbs.ac.in}

\author{Abhishek Pandey}
\address{Abhishek Pandey,
	Discipline of Mathematics,
	School of Basic Sciences,
	Indian Institute of Technology  Bhubaneswar,
	Argul, Bhubaneswar, PIN-752050, Odisha (State),  India.}
\email{ap57@iitbbs.ac.in}


\thanks{}

\maketitle
\pagestyle{myheadings}
\markboth{Vasudevarao  Allu and Abhishek Pandey}{On the Generalized Zalcman Conjecture}

\thanks{}

\maketitle

\begin{abstract}
Let $\mathcal{S}$ denote the class of analytic and univalent ({\it i.e.}, one-to-one) functions $	f(z)= z+\sum_{n=2}^{\infty}a_n z^n$ in the unit disk $\mathbb{D}=\{z\in \mathbb{C}:|z|<1\}$. For $f\in \mathcal{S}$,  In 1999, Ma proposed the generalized Zalcman conjecture that  $$|a_{n}a_{m}-a_{n+m-1}|\le (n-1)(m-1),\,\,\,\mbox{ for } n\ge2,\, m\ge 2,$$
with equality  only for the Koebe function $k(z) = z/(1 - z)^2$ and its rotations. In the same paper, Ma \cite{Ma-1999}  asked  for what positive real values of $\lambda$ does the following inequality hold?
\begin{equation}\label{conjecture}
|\lambda a_na_m-a_{n+m-1}|\le \lambda nm -n-m+1 \,\,\,\,\, (n\ge 2, \,m\ge3).
\end{equation}
Clearly equality holds for  the Koebe function $k(z) = z/(1 - z)^2$ and its rotations. In this paper, we  prove the inequality (\ref{conjecture}) for $\lambda=3, n=2, m=3$. Further, we provide a geometric condition on extremal function maximizing (\ref{conjecture}) for $\lambda=2,n=2, m=3$.
	\end{abstract}

	\section{Introduction}
	Let $\mathcal{H}$ denote the class of analytic functions in the unit disk $\mathbb{D}:=\{z\in\mathbb{C}:\, |z|<1\}$.  Let $\mathcal{A}$ be the class of functions 
	$f\in \mathcal{H}$ such that $f(0)=0$ and $f'(0)=1$, and denote by $\mathcal{S}$ 
	 the class of functions $f\in\mathcal{A}$ which are univalent ({\it i.e.}, one-to-one) in $\mathbb{D}$. 
	Thus, $f\in\mathcal{S}$ has the following representation
	\begin{equation}\label{p4_i001}
		f(z)= z+\sum_{n=2}^{\infty}a_n z^n,\, |z|<1.
	\end{equation}

	In the late 1960's, Zalcman posed the conjecture that if  $f\in\mathcal{S}$, then
	\begin{equation}\label{p4_i005}
		|a_n^2-a_{2n-1}|\le (n-1)^2 \quad\mbox{ for } n\ge2,
	\end{equation}
	with equality  only for the Koebe function $k(z)=z/(1-z)^2$, or its rotation.
	It is important to note that the  Zalcman conjecture implies the celebrated Bieberbach conjecture $|a_n|\le n$ for $f\in\mathcal{S}$  
	(see \cite{Brown-Tsao-1986}), and a  well-known consequence of the area theorem  shows that  (\ref{p4_i005}) holds  for 
	$n=2$ (see \cite{Duren-book-1983}). The Zalcman conjecture remains an  open problem, even after de Branges' proof of the Bieberbach conjecture \cite{Branges-1985}.\\
	
	 For $f\in\mathcal{S}$, Krushkal has proved the Zalcman conjecture for $n=3$ (see \cite{Krushkal-1995}), and recently 
	for $n=4,5$ and $6$ (see \cite{Krushkal-2010}). For a simple and elegant proof of the Zalcman conjecture for the case $n=3$, see \cite{Krushkal-2010}. However, the Zalcman conjecture for $f\in \mathcal{S}$ is still  open for $n>6$. On the other hand,  using complex geometry and universal Teichm\"{u}ller spaces, Krushkal   claimed  in an unpublished work \cite{Krushkal-Unpublished} to have proved the Zalcman conjecture  for all $n\ge 2$. Personal discussions with Prof. Krushkal indicates that there is a gap in the proof of  Krushkal's unpublished work \cite{Krushkal-Unpublished}, and so the Zalcman conjecture  remains open for the class $\mathcal{S}$ for $n>6$.\\ 

	 If $f\in \mathcal{S}$, then the coefficients of $[f(z^2)]^{1/2}$ and $1/f(1/z)$ are polynomials in $a_n$, which contains the expression of the form $\lambda a_n^2-a_{2n-1}$,  pointed out by Pfluger \cite{Pfluger-1976}.

\subsection{Generalized Zalcman conjecture}	In 1999, Ma \cite{Ma-1999} proposed the following generalized Zalcman conjecture: If $f\in \mathcal{S}$, then $$|a_{n}a_{m}-a_{n+m-1}|\le (n-1)(m-1),\,\,\, n\ge 2, m\ge 2$$
	Clearly, for $n=m$, the generalized Zalcman conjecture reduces to the Zalcman conjecture.  which remains an open problem till date. However Ma  \cite{Ma-1999}  proved this generalized Zalcman conjecture for classes $\mathcal{S}^{*}$ and $\mathcal{S}_{\mathbb{R}}$, where   $\mathcal{S}_{\mathbb{R}}$  denotes 
	the class of 
	all functions in $\mathcal{S}$ with real coefficients.\\
	
Further, Ma \cite{Ma-1999}  asked  for what positive real values of $\lambda$ does the following inequality hold?
	\begin{equation}\label{conjecture-2}
		|\lambda a_na_m-a_{n+m-1}|\le \lambda nm -n-m+1 \,\,\,\,\, (n,\,m=2,3,\ldots).
	\end{equation}
Clearly equality holds for  the Koebe function $k(z) = z/(1 - z)^2$ and its rotations.

\begin{rem}\label{Remark}
It is obvious that if (\ref{conjecture-2}) holds, then $\lambda nm -n-m+1\ge0$ {\it i.e.}, 
	$$\lambda\ge\frac{n+m-1}{nm}.$$
Note that if $\lambda\ge\dfrac{n+m-1}{nm}$, $\mu\ge\lambda$ and (\ref{conjecture-2}) holds, then
	\begin{eqnarray*}
		|\mu a_na_m-a_{n+m-1}|&\le&(\mu-\lambda)|a_na_m|+|\lambda a_na_m-a_{n+m-1}|\\ \nonumber&\le&(\mu-\lambda)nm+\lambda nm-n-m+1\,\,\,\,\, (\mbox{since } |a_n|\le n)\\ \nonumber&=&\mu nm-n-m+1.
	\end{eqnarray*}
	That is if \eqref{conjecture-2} holds for some $\lambda$, then it holds for every $\mu\ge \lambda$.
\end{rem}
In view of Remark \ref{Remark}, it is natural to consider the following problem
	\begin{prob}
	Let $f\in \mathcal{S}$, fix $n,m$ and consider the following set
	$\{\lambda: \lambda \mbox{ satisfies }\eqref{conjecture-2}\}$. The problem is to find the infimum of this set.
	\end{prob} 
 At this juncture we must mention that the Zalcman conjecture and its other generalized form has been proved for some subclasses of $\mathcal{S}$, such as starlike functions, typically real functions, close-to-convex functions \cite{Brown-Tsao-1986,Ma-1988,Ponnusamy-2017}. For basic properties of starlike functions, typically real functions and close-to-convex functions we refer to \cite{Duren-book-1983}. 

\section{Main results}
	Using the variational method, 
together with the Bombieri method \cite{Bombieri-1968}, 
we prove the following two results.
\medskip

\begin{thm}\label{Thm-1}
	
	Let $f\in \mathcal{S}$ be an extremal function for the extremal problem $|2a_2a_3-a_4|$ and let $\Gamma$ be the image of $|z|=1$ under $1/f(z)$. If $Re\, a_2>0$ and $Im\, a_2\ne 0$ then $\Gamma$ lies either in upper or lower half plane.
\end{thm}
\medskip

\begin{thm}\label{Thm-2}
	
	Let	$f\in \mathcal{S}$ be given by $f(z)=z+\sum_{n=2}^{\infty} a_n z^n$ then 
	\begin{equation}\label{eq-p4-1-b}
		|3a_2a_3-a_4|\le 14,
	\end{equation}
	with equality  only for  functions of the form
	$$\frac{z}{(1-e^{i \theta}z)^2},\,\,\, \mbox{ where }\,\, \theta \mbox{ is real}.$$
\end{thm}

For the proof of Theorem \ref{Thm-2}, we follow the technique of Ozawa \cite{Ozawa-1964}.

	\section{Preliminaries} In this section we discuss some preliminary ideas which will be useful to prove our main results.
	\subsection{Schiffer's Variational method}
	In 1938, Schiffer developed variational method as a tool for treating the extremal problems arising in univalent function theory. Schiffer's boundary variation \cite{Schiffer-1938,Schiffer-1938-2} which is applicable to very general extremal problems, showed that any function in $\mathcal{S}$ which maximizes
	$|a_{n}|$ must map the disk onto the complement of a single analytic arc \cite{Schiffer-1968} which
	lies on the trajectory of a certain quadratic differential. The omitted arc was found
	to have monotonic modulus and other nice properties. Bieberbach's conjecture
	asserted that this must be a radial half-line. In 1955, Garabedian and Schiffer
	\cite{Schiffer-1955} finally succeeded in using this approach, in combination with Loewner's
	method, to prove that $|a_{4}|\le4$. The work in \cite{Schiffer-1955} gives a general method to
	attack  coefficient problems for univalent functions, but involves a great amount of computational work. Later in 1960, Z. Charzynski and M. Schiffer \cite{Schiffer-1960} gave a greatly simplified proof of
	$|a_4|\le 4$. In \cite{Schiffer-1960}, the authors give a new proof of the Grunsky inequality based on variational methods. \\
	
	 In the application of Schiffer's variational method, the complement $\Gamma$ of the range of an extremal function consists of analytic curves satisfying a differential inequality of the form $Q(w)dw^2<0$.  More precisely, the Schiffer's differential equation is of the  form 
	$$\left(\frac{dw}{dt}\right)^2\,\frac{Q(w)}{w^{n+1}}=-1.$$ 
	In general, such an expression is called a quadratic differential, where $Q$ is meromorphic and the arcs for which  $Q(w)dw^2>0$ are called its trajectories. In many important cases the function $Q$ is a rational function. The zeros and poles of $Q$ are referred to as the singularities of the quadratic differential. In order to emphasise the geometric point of view, consider a metric $ds^2=|Q(w)||dw|^2$, which is euclidean except at the singularities. The trajectories are the geodesics of this metric. For the detailed study of this quadratic differential and its local and global trajectories we refer to Jenkins \cite{Jenkins-1965} and Strebel \cite{Kurt-1984}.\\
	
	 By Schiffer's variational method, the extremal function of certain extremal problems satisfy a differential equation of the  form
	$$z^2 \frac{w'^2}{w^{m+1}}P^*(w)=Q^*(z),$$
	where $P^*$ is a polynomial in $w$ and $Q^*$ is a rational function in $z$. The exact formulation for this differential equation was established by Schaffer and Spencer \cite{Schaeffer-1943,Schaeffer-1945,Schaeffer-1946,Schaeffer-1949,Schaeffer-1949*} in a series of papers on coefficient regions for univalent functions. The problem they considered is to characterize the sequences $\{a_n\}$ which define such functions and can  be solved if one
	can determine the region $V_n$ in $(2n-2)$-dimensional space to which
	the point $(a_2, a_3,\cdots , a_n)$ is confined. The most likely way to success
	is to determine  the boundary of $V_n$ through the extremal
	properties of the corresponding functions, which can be found by developing  a specific variational method for geometric and rigorous meaning of the interior variation, see \cite{Ahlfors}.

	\medskip
	
	  If $f(z)=w$ maximizes an extremal problem $\mathcal{J}(a_2, a_3, a_4, \ldots, a_n)$, then $f$ satisfies the following differential equation (see \cite{Schaeffer-1946})
	
	$$(zw')^2\sum_{v=2}^{n}A_vw^{-v-1}=B+\sum_{v=1}^{n-1}(B_vz^{-n+v}+\overline{B_v}z^{n-v}),$$
	
	where
	
	$$A_v=\sum_{k=v}^{n}a_k^{(v)}\mathcal{J}_k,\,\,\,\, B_v=\sum_{k=1}^{v}ka_k\mathcal{J}_{n+k-v},\,\,\,\, B=\sum_{k=2}^{n}(k-1)a_k\mathcal{J}_k.$$
	
	  Here $a_k^{(v)}$ are the coefficients of 
	$$f(z)^v=\sum_{k=v}^{\infty}a_k^{(v)}z^k,$$
	
	and 
	
	$$\mathcal{J}_v=\frac{\partial \mathcal{J}}{\partial a_v}=\frac{1}{2}\left(\frac{\partial \mathcal{J}}{\partial x_v}-i\frac{\partial \mathcal{J}}{\partial y_v}\right),\,\,\,\,\,\,\, \mbox{ where } a_v=x_v+iy_v.$$
	
\subsection{Bombieri Method} Bombieri \cite{Bombieri-1968}  has proved a general result about critical
trajectories of a quadratic differential $Q(\xi)d\xi^2$ on the $\xi-$ sphere, arising
from the following problem.
\medskip

Let there be given a quadratic differential $Q (\xi) d\xi^2$
on the $\xi$- sphere. Then a "good" subset $T_0$ ($T_0$ is said to be good if it satisfies a certain connectedness condition) of the set $\overline{T}$  where $T$ is the set of   critical trajectories
of $Q (\xi) d\xi^2$, is a continuously differentiable Jordan arc $J$ on the $\xi-$ sphere. Now the question is can we assert that $J\cap T_0$ is either empty set or a single point
under aferesaid conditions on $J$ ?  The answer to this question is given by Bombieri \cite{Bombieri-1968}.

\begin{Thm}\cite{Bombieri-1968}
	 Let $R$ be the $\xi-$ sphere, $Q (\xi)d\xi^2$ be a quadratic differential on it with at most three distint poles, only one of which has order at least $2$. This point is called $B$. $T_0$ be a connected component of $\overline{T}\setminus B$, and let $J$ be a continuously differentiable Jordan arc on $R$ not containing poles of $Q (\xi) d\xi^2$,  such
	that $B\notin \overline{J}$ and
	$$\Im (Q(\xi)d\xi^2)\ne0\,\,\,\, \mbox{ on } J.$$
	Then $\overline{J}$ can meet $T_0$ at most in one point.
\end{Thm}

\begin{Cor}
	Theorem $A$ remains true if $J$ contains one simple pole $A$ of $Q(\xi)d\xi^2$, provided $T_0$ is the connected component of $\overline{T}\setminus B$ containing $A$, and we have
	$$T_0 \cap \overline{J}=A.$$
\end{Cor}

\begin{rem}
	Theorem $A$ and and its corollary remains true if the condition
	$$\Im (Q(\xi)d\xi^2)\ne0\,\,\,\, \mbox{ on } J.$$
	is weakened to 
	$$\Re(Q(\xi)d\xi^2)<0$$
	at every point where $\Im (Q(\xi)d\xi^2)=0\,\,\,\, \mbox{ on } J.$
\end{rem}
	We will use this in our proof to obtain the differential equation satisfied by the extremal function. The main aim of this paper is to solve generalized Zalcman conjcture for the initial coefficients for functions in $\mathcal{S}$.

\section{ Proof of the main results}

We first  prepare some material which will be used in the proof of our main results.\\

	 Let $a_2=x_2+iy_2$, $a_3=x_3+iy_3$, $a_4=x_4+iy_4$, and  $\mathcal{J}=a_4-\lambda a_2a_3$. Then a simple computation gives
	
	$$\mathcal{J}=(x_4-\lambda x_2x_3+\lambda y_2y_3)+i(y_4-\lambda x_2y_3-\lambda y_2x_3).$$	
	
	Further let
	$$\mathcal{J}_{2}:=\frac{1}{2}\left(\frac{\partial \mathcal{J}}{\partial x_2}-i\frac{\partial \mathcal{J}}{\partial y_2}\right)=\frac{1}{2}(-2\lambda x_3-2i\lambda y_3)=-\lambda a_3$$

	\begin{equation}\label{eq-2}
	\mathcal{J}_{3}:=\frac{1}{2}\left(\frac{\partial \mathcal{J}}{\partial x_3}-i\frac{\partial \mathcal{J}}{\partial y_3}\right)=\frac{1}{2}(-2\lambda x_2-2i\lambda y_2)=-\lambda a_2.
	\end{equation}
	
	\begin{equation}\label{eq-3}
		\mathcal{J}_{4}:=\frac{1}{2}\left(\frac{\partial \mathcal{J}}{\partial x_4}-i\frac{\partial \mathcal{J}}{\partial y_4}\right)=\frac{1}{2}(1+1)=1
	\end{equation}

	\begin{eqnarray}\label{eq-4}
		B&=&a_2\mathcal{J}_2+2a_{3}\mathcal{J}_3+3a_4\mathcal{J}_4\\ \nonumber
		&=&-a_2\lambda a_3-2a_3\lambda a_2+3a_4 \\ \nonumber
		&=&3(a_4-\lambda a_2a_3)
	\end{eqnarray}
	
	\begin{eqnarray}\label{eq-5}
		A_2&=&a_2^{(2)}\mathcal{J}_2
		+a_3^{(2)}\mathcal{J}_3+a_4^{(2)}\mathcal{J}_4\\ \nonumber
		&=&\mathcal{J}_2+2a_2\mathcal{J}_3+(a_2^2+2a_3)\mathcal{J}_4\\ \nonumber
		&=&-\lambda a_3-2\lambda a_2^2+a_2^2+2a_3\\ \nonumber
		&=& (1-2\lambda)a_2^2+(2-\lambda)a_3
	\end{eqnarray}
	
	\begin{eqnarray}\label{eq-6}
		A_3&=&a_3^{(3)}\mathcal{J}_3+a_4^{(3)}\mathcal{J}_4\\ \nonumber
		&=&-\lambda a_2+3a_2\\ \nonumber
		&=& (3-\lambda)a_2
	\end{eqnarray}
	
	\begin{eqnarray}\label{eq-7}
		A_4&=&a_4^{(4)}\mathcal{J}_4=1
	\end{eqnarray}
	
	\begin{eqnarray}\label{eq-8}
		B_1&=&\mathcal{J}_4=1
	\end{eqnarray}
	
	\begin{eqnarray}\label{eq-9}
		B_2&=&\mathcal{J}_3+2a_2\mathcal{J}_4\\ \nonumber
		&=& -\lambda a_2+2a_2\\ \nonumber
		&=&(2-\lambda)a_2
	\end{eqnarray}
	
	\begin{eqnarray}\label{eq-10}
		B_3&=&\mathcal{J}_2+2a_2\mathcal{J}_3+3a_3\mathcal{J}_4\\ \nonumber
		&=&- \lambda a_3-2\lambda a_2^2+3a_3\\ \nonumber
		&=& (3-\lambda)a_3-2\lambda a_2^2.
	\end{eqnarray}

	 Since the functional $a_4-\lambda a_2a_3$ is rotationally invariant, we can consider the extremal problem
	$$\max_{\mathcal{S}}\,\Re \{a_4-\lambda a_2a_3\}.$$
	By Schiffer's variational method, the extremal function satisfies the following differential equation

	\begin{eqnarray}\label{eq-11}
		z^2 \frac{w'^2}{w^5}\left(((2-\lambda)a_3+(1-2\lambda)a_2^2)w^2+(3-\lambda)a_2w+1\right)\\ \nonumber
	=\frac{1}{z^3}\left(1+Pz+Qz^2+Rz^3+Sz^4+Tz^5+z^6\right),
	\end{eqnarray}
	where
	
	\begin{eqnarray}\label{equ-12}
		P&=& (2-\lambda)a_{2},\\[2 mm] \nonumber
		Q&=& (3-\lambda)a_3-2\lambda a_{2}^2,\\[2 mm] \nonumber
		R&=& 3(a_4-\lambda a_2a_3),\\[2 mm] \nonumber
		S&=& \overline{(3-\lambda)a_3-2\lambda a_{2}^2}, \\[2 mm] \nonumber
		T&=& (2-\lambda)\overline{a_2}. 
	\end{eqnarray} 
	
	 We are now ready to give the proofs of Theorem \ref{Thm-1} and Theorem \ref{Thm-2}.\\

\begin{pf}[\textbf{Proof of Theorem \ref{Thm-1}}]
	The proof of this theorem requires Bombieri's method \cite{Bombieri-1968} together with Schiffer's variational method. By Schiffer's variational method the image of $|z|=1$ by any extremal function satisfies
	$$\left(\frac{dw}{dt}\right)^2\frac{1}{w^5}\left(-3a_{2}^2w^2+a_2w+1\right)+1=0,$$
	with a suitable parameter $t$. Take  $Q^*(w)dw^2$ as the associated quadratic differential, so that
	$$Q^*(w)dw^2=-\frac{dw^2}{w^5}\left(-3a_{2}^2w^2+a_2w+1\right).$$
	Let $w=1/\xi$ and $Q(\xi)d\xi^2$ be $Q^*(1/\xi)d(1/\xi)^2$, then
	
	$$Q(\xi)d\xi^2=-\left(-3a_2^2+a_2\xi+\xi^2\right)\frac{d\xi^2}{\xi}.$$ 
	Let $a_2=x_2+iy_2$, and $\xi$ be real. Then we have
	$$\Im Q(\xi)d\xi^2=(6x_2y_2-y_2\xi)\frac{d\xi^2}{\xi}.$$ 
	Since $y_2\ne 0$, $\Im Q(\xi)d\xi^2=0$ only if $\xi=6x_2$,
	In view of Bombieri's Theorem 1, and its Remark, we have
	$$\Re Q(\xi)d\xi^2=(-3y_2^2+3x_2^2-x_2\xi-\xi^2)\frac{d\xi^2}{\xi},$$
	and at $\xi=6x_2$ we have
	$$\Re Q(\xi)d\xi^2=-\frac{1}{6x_2}(3y_2^2+39x_2^2)<0,\,\,\,\,\,(\mbox{ since } x_2>0).$$

	Thus from Bombieri's Theorem 1 and its Remark the image $\Gamma$ of $|z|=1$ by $\xi$ intersects the real axis only at the origin. It is easy to observe that $a_2$ can't be zero for any extremal function. Thus,   $\xi=0$ is a simple pole of $Q(\xi)d\xi^2$. Hence $\xi=0$ is an end point of $\Gamma$. Hence $\Gamma$ must lie in either the upper  or lower half plane.
\end{pf}

\begin{pf}[\textbf{Proof of Theorem \ref{Thm-2}}]
Taking $\lambda=3$ in (\ref{eq-11}), we get that the extremal function satisfies the following differential equation

$$z^2 \frac{w'^2}{w^5}\left((-a_3-5a_{2}^2)w^2+1\right)=g(z)=\frac{1}{z^3}\left(1+Pz+Qz^2+Rz^3+Sz^4+Tz^5+z^6\right), $$

where

\begin{eqnarray}\label{9.1}
	P&=& -a_{2},\\[2 mm] \nonumber
	Q&=& -6a_{2}^2,\\[2 mm] \nonumber
	R&=& 3(a_4-3a_2a_3),\\[2 mm] \nonumber
	S&=& -6\overline{a_{2}^2}, \\[2 mm] \nonumber
	T&=& -\overline{a_{2}}. 
\end{eqnarray} 

The image of the unit circle $|z|=1$ under $w=f(z)$ has at least one finite end point, and a point on $|z|=1$, which corresponds to a finite end point is a double zero of $g(z)$. Therefore $g(z)$ can be rewritten as  

\begin{equation}\label{9.2}
	g(z)=\frac{1}{z^3}(z-E)^2(z^4+Dz^3+Cz^2+Bz+A).
\end{equation}
Comparing  (\ref{9.1}) and (\ref{9.2}), we obtain

\begin{eqnarray*}
	BE^2-2AE&=& -a_{2},\\[2 mm] \nonumber
	A-2EB+E^2C&=& -6a_{2}^2,\\[2 mm] \nonumber
	B-2EC+E^2D&=& 3(a_4-3a_2a_3),\\[2 mm] \nonumber
	C-2ED+E^2&=& -6\overline{a_{2}^2},\\[2 mm] \nonumber
	D-2E&=& -\overline{a_{2}}. 
\end{eqnarray*} 
Further $g$ satisfies the following functional equation
$$g(z)=\overline{g\left(\frac{1}{\overline{z}}\right)}.$$

It is easy to see that 

$$AE^2=\overline{A}\overline{E}^2=1,\,\,\,\, |A|=|E|=1,$$
\begin{eqnarray*}	
	BE^2-2AE&=& \overline{D}-2\overline{E},\\[2 mm]
	A-2EB+E^2C&=& \overline{C}-2\overline{E}\overline{D}+\overline{E}^2,\\[2 mm] 
	B-2EC+E^2D&=& \overline{B}-2\overline{E}\overline{C}+\overline{E}^2\overline{D}.
\end{eqnarray*} 
Since $$BE^2-2AE= \overline{D}-2\overline{E},$$ 
we have,
$$D=\overline{B}\overline{E}^2-2\overline{A}\overline{E}+2E=\overline{B}\frac{1}{\overline{A}}+2E-2\frac{1}{\overline{E}^2}\overline{E}=\overline{B}A+2E-2E,$$
and hence 
\begin{equation}\label{A}
	D=\overline{B}A.
\end{equation}

Further,

$$A-2EB+E^2C= \overline{C}-2\overline{E}\overline{D}+\overline{E}^2,$$ 

which implies that
\begin{eqnarray*}
	C&=&\frac{\overline{C}}{E^2}-\frac{2\overline{E}\overline{D}}{E^2}+\frac{\overline{E}^2}{E^2}-\frac{A}{E^2}+\frac{2EB}{E^2}\\[2 mm]
	&=&\overline{C}A-\frac{2\overline{E}B\overline{A}}{E^2}+\frac{1}{E^4}-\frac{1}{E^4}+\frac{2B}{E}\\[2 mm]
	&=&\overline{C}A-\frac{2B}{A}+\frac{2B}{A}.
\end{eqnarray*}
Therefore,
\begin{equation}\label{B}
	C=\overline{C}A.
\end{equation}
By  summarising  we have the following relations:

$$AE^2=1,\,\,\,\, |A|=|E|=1,\,\,\,\, D=\overline{B}A,\,\,\,\, C=\overline{C}A$$

\begin{eqnarray}\label{1}
	BE^2-2AE&=& -a_{2}\\[2 mm]\label{2}
	A-2EB+E^2C&=& -6a_{2}^2 \\[2 mm]\label{3} 
	B-2EC+E^2D&=& 3(a_4-3a_2a_3). 
\end{eqnarray} 
\medskip
Let $ E=e^{i\theta}\, \mbox{ then } A= e^{-2i\theta},\,\,\,\, C=re^{i\beta},\,\,\,\,  B=se^{i\alpha}.$ From (\ref{B}) we obtain
$e^{-2i\beta}=e^{2i\theta}$, which implies that $e^{i(\beta+\theta)}=e^{p\pi i}$, where $  p\in \mathbb{Z}.$
\medskip
From (\ref{1}), (\ref{2}) we obtain the following relations
\begin{equation}\label{equ-p4-2-1}
	se^{i\alpha}-2e^{-3i\theta}=-a_2 e^{-2i\theta}
\end{equation}
\begin{equation}\label{equ-p4-2-2}
	re^{i(\theta+\beta)}-2se^{i\alpha}+e^{-3i\theta}= -6\,a_2^{2}\,e^{-i\theta}.
\end{equation}
Also using (\ref{A}) in (\ref{3}) we obtain
$$3(a_4-3a_2a_3)=B+\overline{B}-2EC=2\operatorname{Re}(B)-2re^{i(\theta+\beta)}=2\,s\cos\alpha-2\,r\cos p\pi,$$
Thus,
\begin{equation}
(3a_2a_3-a_4)=\frac{2}{3}(r\,\cos\,p\pi-s\,\cos\,\alpha).
\end{equation}
Let $-a_2=|-a_2|e^{i\phi}=|a_2|e^{i\phi}$, then from (\ref{equ-p4-2-1}) we have
\begin{equation}\label{equ-p4-2-3}
	s\,\cos\,\alpha-\cos\,3\theta=|a_2|\,\cos\,(\phi-2\theta)+\cos\,3\theta.
\end{equation}
Also from  (\ref{equ-p4-2-2}) we obtain
\begin{equation}\label{equ-p4-2-4}
	r\,\cos\,p\pi-s\,\cos\,\alpha-s\,\cos\,\alpha+\cos\,3\theta=-6|a_2|^2\,\cos(2\phi-\theta),
\end{equation}
and using (\ref{equ-p4-2-3}) in (\ref{equ-p4-2-4}), we obtain

\begin{equation}\label{equ-p4-2-5}
	r\,\cos\,p\pi-s\,\cos\,\alpha=-6|a_2|^2\,\cos(2\phi-\theta)+|a_2|\,\cos\,(\phi-2\theta)+\cos\,3\theta.
\end{equation}
Let $|a_2|=R$, $0<R\le2$ and consider the function 

$$G(R,\theta,\phi)=-6R^2\,\cos\,(2\phi-\theta)+R\,\cos(2\theta-\phi)+\cos\,3\theta.$$
To find the maximum of $G$, first we need to find critical points. It is easy to see that
$$\frac{\partial G}{\partial R}=-12R\,\cos\,(2\phi-\theta)+\cos\,(2\theta-\phi)=0,$$
$$\frac{\partial G}{\partial \theta}=-6R^2\,\sin\,(2\phi-\theta)+2R\sin\,(\phi-2\theta)-3\sin\, 3\theta=0,$$
and
$$\frac{\partial G}{\partial \phi}=12R^2\,\sin\,(2\phi-\theta)-R\,\sin(\phi-2\theta)=0.$$
Hence

$$-12R\,e^{i(\theta-2\phi)}+e^{i(2\theta-\phi)}=0,$$
This gives us that at the points which gives the maximum of $G(R,\theta,\phi)$, we have $$R=\frac{1}{12}, \mbox{ and }  e^{i\phi}=e^{-i\theta}.$$ 
So at $R=1/12,$ we have

$$G=(-6R^2+R+1)\,\cos\,3\phi\le 1.$$
Now we check the value of $G$ under the condition $e^{i\phi}=e^{-i\theta}$ and obtain

$$G=(-6R^2+R+1)\,\cos\,3\phi=(6R^2-R-1)\, (-\cos\,3\phi)\le 21.$$
So the maximum is attained at $R=2$, which gives $G\le21,$ and so $\Re(3a_2a_3-a_4)\le 14,$ {\it i.e.}, $|3a_2a_3-a_4|\le 14,$ with equality  only for the Koebe function $k(z)=z/(1-z)^2$ and its rotations, which  completes the proof.
\end{pf}

\vspace{4mm}
\subsection*{Acknowledgement} Authors thanks Prof. Hiroshi Yanagihara and Prof. D. K. Thomas  for fruitful discussions and giving  constructive suggestions for  improvements to this paper. The first author thanks SERB-CRG, and the second author thanks Prime Minister's Research Fellowship (Id: 1200297) for their support.

\end{document}